\documentclass[12pt]{amsart}
\usepackage{amssymb}

\expandafter\ifx\csname delta.sty\endcsname\relax \else\endinput\fi
\expandafter\edef\csname delta.sty\endcsname{%
 \catcode`\noexpand\@=\the\catcode`\@\space}

\let\atbefore @

\catcode`\@=11

\newif\ifMag
\let\@ft@\expandafter
\numberwithin{equation}{section}

\newif\ifComments

\def\tod@y{\ifcase\month\or
 January\or February\or March\or April\or May\or June\or July\or
 August\or September\or October\or November\or December\fi\space\,
\number\day,\space\,\number\year}

\def\h@@r{hh}\def\m@n@te{mm}
\def\wh@tt@me{\count@\time\divide\count@ 60\edef\h@@r{\number\count@}%
 \multiply\count@ -60\advance\count@\time\edef
 \m@n@te{\ifnum\count@<10 0\fi\number\count@}}
\def\t@me{\h@@r\/{\rm:}\m@n@te} \let\whattime\wh@tt@me
\let\Today\tod@y \let\nowtime\t@me

\def\ftext#1{{\let\thefootnote\relax\footnotetext{\vsk-.8>\nt #1}}}
\def\em#1{{\itshape #1\/}}

 {\end{trivlist}\egroup\global\@ignoretrue\ifComments\unvbox\commentbox\fi
 \noindent}

\def\gadv{\global\adv} \def\gad#1{\gadv#1\@ne} \def\gadneg#1{\gadv#1-\@ne}
\def\textindent#1{\indent\llap{#1\enspace}\ignorespaces}

\newcount\itemlet
\def\newbi{\itemlet 96} \newbi
\def\bitem{\gad\itemlet\endgraf\hangindent1.5\parindent
 \hglue-.5\parindent\textindent{\upshape\rlap{\char\the\itemlet}\hp{b})}}

\newcount\itemrm

\def\iitem{\gad\itemrm\endgraf\hangindent1.5\parindent\hglue-.5\parindent
 \textindent{\upshape\hp{v}\llap{\romannumeral\the\itemrm})}}

\let\Disp\[ \let\endD\] \let\{\protect

\def\Tag#1{\label{e:#1}\let\notag\relax} 
\def\sh@nd#1#2{\begin{#1*}#2\end{#1*}}
\def\n@t@gs#1#2#3{\let\n@@@l\\ \begin{#1}#2\global\let\d@bl\\
 \gdef\\{\notag\d@bl}#3\notag\global\let\\\d@bl\end{#1}\let\\\n@@@l}

\def\Gather#1\endG{\n@t@gs{gather}{}{#1}}
\def\gAther#1\endG{\sh@nd{gather}{#1}}
\def\Align#1\endA{\n@t@gs{align}{}{#1}}
\def\aLign#1\endA{\sh@nd{align}{#1}}
\def\Alignat#1#2\endAt{\n@t@gs{alignat}{#1}{#2}}
\def\aLignat#1\endAt{\sh@nd{alignat}{#1}}

\def\(#1){\textup{(\ref{e:#1})}}
\def\[{\@ifnextchar:\c@t@sect\c@t@}
\def\c@t@sect:#1]{\ref{s:#1}} \def\c@t@#1]{\ref{t:#1}}

\def\qed{\hbox{}\nobreak\hfill\nobreak{\m@th$\,\square$}}
\def\sk@@p#1{\par\skip@#1\relax\ifdim\lastskip<\skip@\relax\removelastskip
 \vskip\skip@\fi}
\def\demo#1{\sk@@p\medskipamount\nt{\ignore\it #1\unskip.}\enspace
 \ignore}
\def\enddemo{\sk@@p\medskipamount}
 
\def\Pf#1.{\demo{Proof #1}}

\let\bls\baselineskip \let\ignore\ignorespaces \let\adv\advance
\def\vsk#1>{\vskip#1\bls}
\def\vv#1>{\vadjust{\vsk#1>}\ignore}
\def\vvn#1>{\vadjust{\nobreak\vsk#1>\nobreak}\ignore}
\def\vvv#1>{\vskip\z@\vsk#1>\nt\ignore}

\def\mathbox#1{\hbox{\m@th$#1$}}

\let\dsize\displaystyle \let\tsize\textstyle
\let\ssize\scriptstyle \let\sss\scriptscriptstyle
 
\let\vp\vphantom \let\hp\hphantom \let\nt\noindent
\def\hline{\hbox to\hsize}
\let\cline\centerline \let\lline\leftline \let\rline\rightline
\def\nn#1>{\noalign{\vskip#1\p@}} \def\NN#1>{\openup#1\p@}
 
\let\Lim\lim \def\lim{\Lim\limits} \let\Sum\sum \def\sum{\Sum\limits}
 
\let\Prod\prod \def\prod{\Prod\limits} \let\Int\int \def\int{\Int\limits}

\def\~{\leavevmode\@ifnextchar~\m@n@s\@md@sh}
\def\m@n@s~{\raise.15ex\mathbox{-}} \def\@md@sh{\raise.13ex\hbox{--}}
\let\procent\% \def\%#1{\ifmmode\mathop{#1}\limits\else\procent#1\fi}
\let\@ml@t\" \def\"#1{\ifmmode ^{(#1)}\else\@ml@t#1\fi}
\let\@c@t@\' \def\'#1{\ifmmode _{(#1)}\else\@c@t@#1\fi}
\let\colon\: \def\:{^{\vp|}} \def\&{.\kern.1em} \def\^#1{\text{\m@th#1}}

\newif\ifNewskips

\def\Newskips{\global\Newskipstrue
 \gdef\>{\relax\ifmmode\mskip.666667\thinmuskip\relax\else\kern.111111em\fi}
 \gdef\}{\relax\ifmmode\mskip-.666667\thinmuskip\relax\else\kern-.111111em\fi}
 \gdef\){\relax\ifmmode\mskip.333333\thinmuskip\relax\else\kern.0555556em\fi}
 \gdef\]{\relax\ifmmode\mskip-.333333\thinmuskip\relax\else\kern-.0555556em\fi}}
\Newskips

 \def\Im{\mathop{\mathrm{Im}\>}}
 
 \def\Sym{\mathop{\mathrm{Sym}\)}}

\def\1{^{-1}} \def\_#1{_{\Rlap{#1}}}
\def\vst#1{{\lower1.9\p@
 \hbox{\m@th$\bigr|_{\raise.5\p@\hbox{\m@th$\ssize#1$}}$}}}
\def\vrp#1:#2>{{\vrule height#1 depth#2 width\z@}}
\def\vru#1>{\vrp#1:\z@>} \def\vrd#1>{\vrp\z@:#1>}
 
\def\sscr#1{\raise.3ex\hbox{\m@th$\sss#1$}} \def\@@PS{\mathbf{OOPS!!!}}

\def\lsym#1{#1\alb\ldots\relax#1\alb}
\def\lc{\lsym,}

\let\texspace\ \def\ {\ifmmode\alb\fi\texspace}

\def\Line#1{\kern-.5\hsize\hline{\m@th$\dsize#1$}\kern-.5\hsize}
\def\Lline#1{\kern-.5\hsize\lline{\m@th$\dsize#1$}\kern-.5\hsize}
\def\Cline#1{\kern-.5\hsize\cline{\m@th$\dsize#1$}\kern-.5\hsize}
\def\Rline#1{\kern-.5\hsize\rline{\m@th$\dsize#1$}\kern-.5\hsize}

\def\Ll@p#1{\llap{\m@th$#1$}} \def\Rl@p#1{\rlap{\m@th$#1$}}
 \def\Cl@p#1{\llap{\m@th$#1$\hss}}
\def\Llap#1{\mathchoice{\Ll@p{\dsize#1}}{\Ll@p{\tsize#1}}{\Ll@p{\ssize#1}}%
 {\Ll@p{\sss#1}}}
\def\Clap#1{\mathchoice{\Cl@p{\dsize#1}}{\Cl@p{\tsize#1}}{\Cl@p{\ssize#1}}%
 {\Cl@p{\sss#1}}}
\def\Rlap#1{\mathchoice{\Rl@p{\dsize#1}}{\Rl@p{\tsize#1}}{\Rl@p{\ssize#1}}%
 {\Rl@p{\sss#1}}}
 
\def\LRtph#1#2{\setbox\z@\hbox{#1}\dimen\z@\wd\z@\hbox{\hbox to\dimen\z@{#2}}}
\def\LRph#1#2{\LRtph{\m@th$#1$}{\m@th$#2$}}

\def\Lto#1{\setbox\z@\hbox{\m@th$\tsize{#1}$}%
 \mathrel{\mathop{\hbox to\wd\z@{\rightarrowfill}}\limits#1}}
\def\Lgets#1{\setbox\z@\hbox{\m@th$\tsize{#1}$}%
 \mathrel{\mathop{\hbox to\wd\z@{\leftarrowfill}}\limits#1}}

\let\alb\allowbreak

 \let\x\times \let\ox\otimes 
 
\let\le\leqslant \let\ge\geqslant
\let\der\partial \let\8\infty \let\*\star

 \def\vert{\ |\ }

\let\lb\lbrace \let\rb\rbrace

\let\Bbb\mathbb

\let\Cal\mathcal
\let\frak\mathfrak

\def\pms{\raise.25ex\mathbox{\ssize\pm}\>}
\def\mps{\raise.25ex\mathbox{\ssize\mp}\>}

\let\gm\gamma \let\Gm\Gamma 
\let\dl\delta  
 \let\eps\varepsilon \let\epsilon\eps

\let\ka\kappa
\let\la\lambda

 \let\phi\varphi

\let\om\omega \let\Om\Omega 

\def\C{\Bbb C}

\def\Z{\Bbb Z}

\def\Zp{\Z_{\ge 0}}

\def\h@ph{\discretionary{}{}{-}} \def\$#1$-{\,\^{$#1$}\h@ph}

\def\difl/{differential} \def\dif/{difference}
\def\cf.{cf.\ \ignore} \def\Cf.{Cf.\ \ignore}
\def\egv/{eigenvector} \def\eva/{eigenvalue} \def\eq/{equation}
\def\lhs/{the left hand side} \def\rhs/{the right hand side}
\def\Lhs/{The left hand side} \def\Rhs/{The right hand side}
\def\gby/{generated by} \def\wrt/{with respect to} \def\st/{such that}
\def\resp/{respectively} \def\off/{offdiagonal} \def\wt/{weight}
\def\pol/{polynomial} \def\rat/{rational} \def\tri/{trigonometric}
\def\fn/{function} \def\var/{variable} \def\raf/{\rat/ \fn/}
\def\inv/{invariant} \def\hol/{holomorphic} \def\hof/{\hol/ \fn/}
\def\mer/{meromorphic} \def\mef/{\mer/ \fn/} \def\mult/{multiplicity}
\def\sym/{symmetric} \def\perm/{permutation}
\def\rep/{representation} \def\irr/{irreducible} \def\irrep/{\irr/ \rep/}
\def\hom/{homomorphism} \def\aut/{automorphism} \def\iso/{isomorphism}
\def\lex/{lexicographical} \def\as/{asymptotic} \def\asex/{\as/ expansion}
\def\ndeg/{nondegenerate} \def\neib/{neighbourhood} \def\deq/{\dif/ \eq/}
\def\hw/{highest \wt/} \def\gv/{generating vector} \def\eqv/{equivalent}
\def\msd/{method of steepest descend} \def\pd/{pairwise distinct}
\def\wlg/{without loss of generality} \def\Wlg/{Without loss of generality}
\def\onedim/{one-dim\-en\-sion\-al} \def\fd/{fi\-ni\-te-dim\-en\-sion\-al}
\def\qcl/{quasiclassical} \def\hwv/{\hw/ vector}
\def\hgeom/{hyper\-geo\-met\-ric} \def\hint/{\hgeom/ integral}
\def\hwm/{\hw/ module} \def\emod/{evaluation module} \def\Vmod/{Verma module}
\def\symg/{\sym/ group} \def\sol/{solution} \def\eval/{evaluation}
\def\anf/{analytic \fn/} \def\anco/{analytic continuation}
\def\qg/{quantum group} \def\qaff/{quantum affine algebra}

\def\Rm/{\^{$R$-}matrix} \def\Rms/{\^{$R$-}matrices}
\def\YB/{Yang-Baxter \eq/}
\def\Ba/{Bethe ansatz} \def\Bv/{Bethe vector} \def\Bae/{\Ba/ \eq/}
\def\KZv/{Knizh\-nik-Zamo\-lod\-chi\-kov} \def\KZvB/{\KZv/-Bernard}
\def\KZ/{{\sl KZ\/}} \def\qKZ/{{\sl qKZ\/}}
\def\KZB/{{\sl KZB\/}} \def\qKZB/{{\sl qKZB\/}}
\def\qKZo/{\qKZ/ operator} \def\qKZc/{\qKZ/ connection}
\def\KZe/{\KZ/ \eq/} \def\qKZe/{\qKZ/ \eq/} \def\qKZBe/{\qKZB/ \eq/}

\def\LPT/{Laboratoire de Physique Th\'eorique ENSLAPP}
\def\ENSLyon/{\'Ecole Normale Sup\'erieure de Lyon}
\def\DMS/{Department of Mathematics, Faculty of Science}
\def\DMO/{\DMS/, Osaka University}
\def\DMOaddr/{Toyonaka, Osaka 560, Japan}
\def\dmoemail/{vt@math.sci.osaka-u.ac.jp}
\def\SPb/{St\&Petersburg}
\def\home/{\SPb/ Branch of Steklov Mathematical Institute}
\def\homeaddr/{Fontanka 27, \SPb/ \,191011, Russia}
\def\homemail/{vt@pdmi.ras.ru}
\def\absence/{On leave of absence from \home/}
\def\UNC/{Department of Mathematics, University of North Carolina}
\def\ChH/{Chapel Hill}
\def\UNCaddr/{\ChH/, NC 27599, USA} \def\avemail/{av@math.unc.edu}
\def\grant/{NSF grant DMS\~9501290}	
\def\Grant/{Supported in part by \grant/}

\def\Aomoto/{K\&Aomoto}
\def\Dri/{V\]\&G\&Drin\-feld}
\def\Fadd/{L\&D\&Fad\-deev}
\def\Feld/{G\&Felder}
\def\Fre/{I\&B\&Fren\-kel}
\def\Gustaf/{R\&A\&Gustafson}
\def\Kazh/{D\&Kazhdan} \def\Kir/{A\&N\&Kiril\-lov}
\def\Kor/{V\]\&E\&Kore\-pin}
\def\Lusz/{G\&Lusztig}
\def\MN/{M\&Naza\-rov}
\def\Resh/{N\&Reshe\-ti\-khin} \def\Reshy/{N\&\]Yu\&Reshe\-ti\-khin}
\def\Skl/{E\&K\&Sklya\-nin}
\def\SchV/{V\]\&\]V\]\&Schecht\-man} \def\Sch/{V\]\&Schecht\-man}
\def\Takh/{L\&A\&Takh\-tajan}
\def\VT/{V\]\&Ta\-ra\-sov} \def\VoT/{V\]\&O\&Ta\-ra\-sov}
\def\Varch/{A\&\]Var\-chenko} \def\Varn/{A\&N\&\]Var\-chenko}

\def\AMS/{Amer.\ Math.\ Society}
\def\CMP/{Comm.\ Math.\ Phys.{}}
\def\DMJ/{Duke.\ Math.\ J.{}}
\def\Inv/{Invent.\ Math.{}} 
\def\IMRN/{Int.\ Math.\ Res.\ Notices}
\def\JPA/{J.\ Phys.\ A{}}
\def\JSM/{J.\ Soviet\ Math.{}}
\def\LMP/{Lett.\ Math.\ Phys.{}}
\def\LMJ/{Leningrad Math.\ J.{}}
\def\LpMJ/{\SPb/ Math.\ J.{}}
\def\SIAM/{SIAM J.\ Math.\ Anal.{}}
\def\SMNS/{Selecta Math., New Series}
\def\TMP/{Theor.\ Math.\ Phys.{}}
\def\ZNS/{Zap.\ nauch.\ semin. LOMI}

\def\ASMP/{Advanced Series in Math.\ Phys.{}}

\def\AMSa/{AMS \publaddr Providence}
\def\Birk/{Birkh\"auser}
\def\CUP/{Cambridge University Press} \def\CUPa/{\CUP/ \publaddr Cambridge}
\def\Spri/{Springer-Verlag} \def\Spria/{\Spri/ \publaddr Berlin}
\def\WS/{World Scientific} \def\WSa/{\WS/ \publaddr Singapore}

\csname delta.sty\endcsname

\newif\ifUS
\ifx\USversion\relax\UStrue\fi

\ifUS
\fi

\evensidemargin\oddsidemargin
\newtheorem{theorem}{Theorem}
\newtheorem{lemma}[theorem]{Lemma}
\newtheorem{corollary}[theorem]{Corollary}

\newenvironment{abstr}{\begingroup\narrower\small
\nt{\sc Abstract.}\enspace\ignorespaces}{\endgraf\endgroup}

\def\fratop{\genfrac{}{}{0pt}1}
\def\Ref#1{{\rm(\ref{#1})}}

\def\Ib{\bar I}

\def\Cc{\Cal C}
\let\Oc O

\def\h{\frak h}
\def\n{\frak n}

\def\gl{\frak{gl}}
\def\glt{\gl_2}

\def\Ih{\hat I}
\def\V-{\>\hbox{$\=V\}$-}}

\title[\smash{Identities for Hypergeometric Integrals of Different Dimensions}]
{Identities for Hypergeometric Integrals\\[3pt] of Different Dimensions}

\author[\smash{V\]\&Tarasov and A\&\]Varchenko}]
{\vbox{}V\]\&Tarasov$^\star$ \>and \;A\&\]Varchenko$^\diamond$}

\begin{document}

\hrule width0pt

\maketitle

\ftext{\mathsurround 0pt
$\]^\star\)$Supported in part by RFFI grant 02\)\~\)01\~\)00085a
\>and \,CRDF grant RM1\~\)2334\)\~MO\)\~\)02\vv.2>\\
\hp{$^*$}{\normalsize\sl E-mail\/{\rm:} \homemail/}\vv.1>\\
${\]^\diamond\)}$Supported in part by NSF grant DMS\)\~\)0244579\vv.2>\\
\vv-1.2>
\hp{$^*$}{\normalsize\sl E-mail\/{\rm:} anv@email.unc.edu}}

\begin{center}
{\it\home/\\[2pt]\homeaddr/\\[6pt]
$^\diamond$\UNC/ at \ChH/\\[2pt]\UNCaddr/}

\vsk1.5>
{\sl May, 2003}
\end{center}

\vsk1.2>

\begin{abstr}
Given complex numbers $m_1\),l_1$ and nonnegative integers $m_2\),l_2$, such
that $m_1+m_2 = l_1+l_2$, we define \$l_2\)$-dimensional hypergeometric
integrals $I_{a,b}(z; m_1, m_2, l_1, l_2)$, $a\),b=0\lc\min\)(m_2\),l_2)$,
depending on a complex parameter $z$. We show that
$I_{a,b}(z;m_1,\alb m_2,l_1,l_2)=I_{a,b}(z;l_1,l_2,m_1,m_2)$,
thus establishing an equality of $l_2$ and \$m_2\)$-dimensional
integrals. This identity allows us to study asymptotics of the integrals with
respect to their dimension in some examples. The identity is based on the
$(\gl_k\>,\gl_n)$ duality for the \KZ/ and dynamical differential equations.
\end{abstr}

\thispagestyle{empty}

\section{Introduction}

Let $\ka$ be a positive number.
Let $m_1, l_1$ be complex numbers and $m_2, l_2$ nonnegative integers
such that
\vvn-.4>
$$
m_1+\>m_2\,=\,l_1+\>l_2\,.
\vv.2>
$$
We say that an integer $a$ is \em{admissible} with respect to $m_2, l_2$ if
\vvn.2>
$$
0\le a\le\min\)(m_2, l_2)\,.
$$

\vsk.1>
For a pair of admissible numbers $a, b$ we define a function
$I_{a,b}(z; m_1, m_2, l_1, l_2)$
of one complex variable $z$. The function is defined as
an \$l_2\)$-dimensional hypergeometric integral of the form
\begin{align}
\label{main integral}
I_{a,b}(z; m_1, m_2, l_1, l_2)\,={}\!\] &
\\[6pt]
{}=\,C_{b}(m_1,m_2,l_1,l_2) &
\int_{\gm_{l_2,\)b}(z)\!}\!\Phi_{l_2}(t, z; m_1, m_2)^{1/\ka}
\,w_{l_2-a,\)a}(t, z)\,dt^{\>l_2}\,.
\notag
\end{align}
Here $t=(t_1,\dots, t_{l_2})$, \,$dt^{\>l_2}=dt_1\wedge\ldots\wedge dt_{l_2}$.
The constant $C_{b}(m_1,m_2, l_1,l_2)$ and the functions
$\Phi_{l_2}(t, z; m_1, m_2)$ and $w_{l_2-a,\)a}(t, z)$ are defined below.
The \$l_2\)$-dimensional integration contour $\gm_{l_2,b}(z)$ lies in
$\C^{\>l_2}$ and is also defined below.

\goodbreak
The \em{master} function $\Phi_l$ is defined by the formula
\begin{equation}
\label{master}
\Phi_l(t_1,\dots , t_l, z; m_1, m_2)\,=\,e^{-\!\sum_{u=1}^l\!t_u}\,
\prod_{u=1}^l\,(-\)t_u)^{-m_1}\>(z-t_u)^{-m_2}
\prod_{1\le u<v\le l}(t_u\]-t_v)^2\,.
\end{equation}
The \em{weight} function $w_{l-a,\)a}$ is defined by the formula
\begin{equation}
\label{weight}
w_{l-a,\)a}(t_1,\dots,t_l, z) \,=\,
\Sym \biggl[\ \prod_{u=1}^{l-a}\,\frac 1{-\)t_u}\
\prod_{u=l-a+1}^l\,\frac 1{z - t_u}\ \biggr]
\end{equation}
where \>$\Sym f(t_1, \dots , t_{l})\>=\sum_{\sigma\in S_{l}}\,
f(t_{\sigma_1}, \dots , t_{\sigma_{l}})$.
We define the integral in \Ref{main integral} by analytic
continuation from the region $\Im z>0$. For $\Im z>0$ the integration contour
$\gm_{l,b}(z)$ is shown in the picture.

\vsk.5>
\vbox{\begin{center}
\begin{picture}(222,85)
\put(40,15){\oval(10,10)[l]}
\put(40,15){\circle*{3}}
\put(40,10){\vector(1,0){152}}
\put(40,20){\line(1,0){152}}
\qbezier[16](46,5)(118,5)(190,5)
\qbezier[16](46,25)(118,25)(190,25)
\qbezier[4](30,15)(30,25)(40,25)
\qbezier[4](40,5)(30,5)(30,15)
\put(40,15){\oval(30,30)[l]}
\put(40,0){\vector(1,0){152}}
\put(40,30){\line(1,0){152}}

\put(31,65){\oval(10,10)[l]}
\put(31,65){\circle*{3}}
\put(31,60){\vector(1,0){161}}
\put(31,70){\line(1,0){161}}
\qbezier[17](37,55)(113,55)(190,55)
\qbezier[17](37,75)(113,75)(190,75)
\qbezier[4](21,65)(21,75)(31,75)
\qbezier[4](31,55)(21,55)(21,65)
\put(31,65){\oval(30,30)[l]}
\put(31,50){\vector(1,0){161}}
\put(31,80){\line(1,0){161}}

\put(15,5){$0$}
\put(197,30){${t_{l}}$}
\put(197,20){${t_{b+1}}$}

\put(6,55){$z$}
\put(197,80){${t_{b}}$}
\put(197,70){${t_{1}}$}
\end{picture}
\vsk>
Picture 1. The contour $\gm_{l,b}(z)$.
\end{center}}

\nt
It has the form
\vvn-.3>
$$
\gm_{l,b}(z)\,=\,\lb\,(t_1,\dots, t_l) \in \C^{\>l}
\ |\ \,t_u\in\Cc_u\,,\ \;u = 1,\dots, l\;\rb\,.
$$
Here $\Cc_u$, $u=1,\dots, l$, are non-intersecting oriented loops in $\C$.
The first $b$ loops start at $+\)\infty$, go around $z$, and return to
$+\)\infty$. For $1\le u<v\le b$, the loop $\Cc_u$ lies inside the loop
$\Cc_v$. The last $l-b$ loops start at $+\)\infty$, go around $0$, and
return to $+\)\infty$. For $b+1\le u<v\le l$, the loop $\Cc_u$ lies inside
the loop $\Cc_v$.
\vsk.2>
If $\Im z>0$, then we fix a univalued branch of
$\Phi_{l}(t, z; m_1, m_2)^{1/\ka}$ over $\gm_{l,b}(z)$ by fixing
the arguments of all factors of $\Phi_l$. Namely, we assume that
at the point of $\gm_{l,b}(z)$ where all numbers
$z-t_1,\dots, z-t_b,\,-\)t_{b+1},\dots, -\)t_{l}$ are positive we have
\begin{enumerate}
\item[]
$\arg\)(-\)t_u)\in (\)-\)\pi, 0)$ for $u = 1,\dots, b$, and
$\arg\)(-\)t_u) = 0$ for $u = b+1, \dots , l$;
\item[]
$\arg\)(z-t_u) = 0$ for $u = 1, \dots , b$, and
$\arg\)(z-t_u) \in (0, \pi)$ for $u = b+1, \dots , l$;
\item[]
$\arg\)(t_u\]-t_v) = 0$ for $u< v\le b$ and for $b<u<v$;
\item[]
$\arg\)(t_u\]-t_v)\in(0,\pi)$ for $u\le b<v$.
\end{enumerate}
Set
\vvn-.8>
\begin{align}
C_b( &m_1,m_2,l_1,l_2)\,=\,\ka^{(l_1+1)l_2/\ka}\,e^{-\pi i\)b\)l_2/\ka}
\,\bigl(2i\,\Gm(-1/\ka)\bigr)^{-l_2}\>\x{}
\notag
\\[10pt]
& \!\]{}\x\;\prod_{j=0}^{b-1}\;\frac1{\sin\bigl(\pi(j+1)/\ka\bigr)}\;
\prod_{j=0}^{l_2-b-1}\frac1{\sin\bigl(\pi(j+1)/\ka\bigr)}\,\,
\prod_{j=0}^{l_2-1}\;
\frac{\Gm\bigl(1+(m_1-j)/\ka\bigr)}{\Gm\bigl(1+(j+1)/\ka\bigr)}\;.
\notag
\end{align}

\vsk.5>
\begin{theorem}
\label{first}
Let $m_1, l_1$ be complex numbers and $m_2, l_2$ nonnegative integers,
such that $m_1 + m_2 = l_1 + l_2$. Let $\ka$ be a generic positive number.
Then for any $a, b$, admissible with respect to $m_2, l_2$, we have
$$
I_{a,b}(z; m_1, m_2, l_1, l_2)\,=\,I_{a,b}(z; l_1, l_2, m_1, m_2)\,.
$$
\end{theorem}

\goodbreak
\vsk.5>

This theorem claims that an \$l_2\)$-dimensional integral equals
an \$m_2\)$-dimensional integral. Note that in the first integral the numbers
$m_1, m_2$ are exponents of its master function, while in the second integral
$m_2$ is its dimension and $m_1$ is not present explicitly.

It is well known that studying asymptotics of integrals with respect
to their dimension is an interesting problem appearing, for instance,
in the theory of orthogonal polynomials and in the matrix models.
The duality of the theorem allows us to study asymptotics of integrals with
respect to their dimension in some examples. Namely, assume that in the
4\)-tuple $(m_1, m_2, l_1, l_2)$ the nonnegative integer $l_2$ tends
to infinity, while the numbers $m_1$ and $m_2$ remain fixed.
Then $I_{a,b}(z;m_1,m_2,l_1,l_2)$ is an integral of growing dimension $l_2$,
whose master function has fixed exponents $m_1, m_2$. At the same time,
$I_{a,b}(z; l_1,l_2, m_1,m_2)$ is an integral of the fixed dimension $m_2$,
whose master function has growing exponents $l_1, l_2$. The asymptotics of
$I_{a,b}(z; l_1,l_2, m_1,m_2)$ can be calculated using the steepest descent
method, see a more precise statement below.

\vsk.2>
To prove Theorem~\ref{first} we show that the matrices:
$$
\bigl(I_{a,b}(z;m_1,m_2,l_1,l_2)\bigr)_{0\le a,\)b\)\le\min\)(m_2,\)l_2)}
\quad\text{and}\quad
\bigl(I_{a,b}(z;l_1,l_2,m_1,m_2)\bigr)_{0\le a,\)b\)\le\min\)(m_2,\)l_2)}
$$
satisfy the same first order linear differential equation with respect to $z$,
and have the same asymptotics as $z$ tends to infinity,
see Lemmas~\ref{lem:eq} and \ref{lem:as}. Hence, the matrices are equal.

The fact that both matrices satisfy the same differential equation is based on
the duality of the \KZ/ and dynamical differential equations for $\gl_k$ and
$\gl_n$ \cite{TV1}\). Namely, according to \cite{FMTV} the \KZ/ and dynamical
differential equations for $\gl_k$ have the hypergeometric solutions, and
the \KZ/ and dynamical differential equations for $\gl_n$ also have
the hypergeometric solutions. The $(\gl_k\>,\gl_n)$ duality allows one to
conclude that both sets of solutions solve the same differential equations.
Thus, the two sets of solutions, in principle, can be identified. Theorem
\ref{first} is a realization of this idea for $k=n=2$. We will discuss
the case of an arbitrary pair $k,n$ in a separate paper.

\vsk.2>
In Section~2 we give a proof of Theorem~\ref{first}. In Section~3
we discuss the $(\glt\>,\glt)$ duality for the \KZ/ and dynamical differential
equations. In Section~4 we construct hypergeometric solutions of the \KZ/
and dynamical differential equations and relate them to the integrals
$I_{a,b}(z;m_1,m_2,l_1,l_2)$. In Section~5 we calculate an example of
asymptotics of hypergeometric integrals with respect to their dimension.

\vsk.2>

The authors thank Y\]\&Markov for numerous useful discussions.

\section{Proof of Theorem \ref{first}}
Introduce matrices $A$ and $B(m_1,m_2,l_1,l_2)$ with entries
\vvn.3>
\begin{gather}
A_{a,b}\,=\,a\,\dl_{a,b}\,,
\notag
\\[6pt]
\begin{aligned}
B_{a,b}(m_1,m_2,l_1,l_2)\,&{}=\,
\bigl(2a^2-a(2l_2+m_2-m_1)+m_2l_2\bigr)\,\dl_{a,b}\,+{}
\\[4pt]
& \>{}+\,a(l_2\]-m_1\]-a)\,\dl_{a-1,b}\,-\,
(m_2-a)(l_2-a)\,\dl_{a+1,b}\,,
\end{aligned}
\notag
\end{gather}
$a,b=0,\dots,\min\)(m_2,l_2)$.

\vsk.5>
\begin{lemma}
\label{lem:eq}
The matrix $\Ih(z;m_1,m_2,l_1,l_2)\>=\>
\bigl(I_{a,b}(z;m_1,m_2,l_1,l_2)\bigr)_{0\le a,\)b\)\le\min\)(m_2,\)l_2)}$
obeys the ordinary differential equation
\vvn.5>
\begin{equation}
\label{eq:Ihat}
\biggl(\,\ka\,\frac{\der}{\der z}\,+\,\frac{B(m_1,m_2,l_1,l_2)}z\,+\,A\,\biggr)
\,\Ih(z;m_1,m_2,l_1,l_2)\,=\,0\,.
\end{equation}
\end{lemma}
\vsk.3>\nt
The lemma admits a straightforward proof using integration by parts.
In Section~4 we give another proof of Lemma~\ref{lem:eq} using the
hypergeometric solutions of the \KZ/ and dynamical differential equations.

\vsk.5>
\begin{lemma}
\label{lem:as}
Fix a positive number $\eps$ less than $\pi/2$.
If $z$ tends to infinity inside the sector \>$\eps<\arg\)z<\pi-\eps$, then
\vvn.4>
\begin{align}
I_{a,b}(z;m_1,m_2,l_1,l_2)\, &{}=\,
e^{-bz/\ka}\,(\ka/z)^{(2b^2-b\)(2l_2+m_2-m_1)+m_2l_2)/\ka}\,
e^{\pi i\)b\)(m_1-\)l_2)/\ka}\,\x{}
\notag
\\[4pt]
& \>{}\x\,\prod_{j=0}^{b-1}\,
\frac{\Gm\bigl(1+(j+1)/\ka\bigr)\,\Gm\bigl(1+(m_1-l_2+j+1)/\ka\bigr)}
{\Gm\bigl(1+(m_2-j)/\ka\bigr)\,\Gm\bigl(1+(l_2-j)/\ka\bigr)}
\ \bigl(\dl_{a,b}+\Oc(z^{-1})\bigr)\,.\!
\notag
\end{align}
\end{lemma}

\vsk.5>
\begin{proof}
For a positive number $\ka$, complex number $m$, and a nonnegative integer $l$,
consider the Selberg type integral
\vvn.3>
$$
J_l(m)\,=\,
\int_{\dl_l}\,e^{-\!\sum_{u=1}^l\!s_u/\ka}\,\prod_{u=1}^l(-\)s_u)^{-1-m/\ka}
\prod_{1\le u<v\le l}(s_u\]-s_v)^{2/\ka}\,ds^k\,.
\vv.3>
$$
The integration contour $\dl_l$ has the form
$\dl_l\>=\>\lb\,(s_1,\dots,s_{l})\in\C^{\>l}\vert s_u\in\Cc_u\,,
\ u=1,\dots, l\,\rb$, see the picture.

\vbox{\begin{center}
\begin{picture}(224,35)
\put(40,15){\oval(10,10)[l]}
\put(40,15){\circle*{3}}
\put(40,10){\vector(1,0){152}}
\put(40,20){\line(1,0){152}}
\qbezier[16](46,5)(118,5)(190,5)
\qbezier[16](46,25)(118,25)(190,25)
\qbezier[5](30,15)(30,25)(40,25)
\qbezier[4](40,5)(30,5)(30,15)
\put(40,15){\oval(30,30)[l]}
\put(40,0){\vector(1,0){152}}
\put(40,30){\line(1,0){152}}
\put(15,5){$0$}
\put(197,30){$s_l$}
\put(197,20){$s_1$}
\end{picture}
\vsk>
Picture 2. The contour $\dl_l$.
\end{center}}
\vsk.3>

\noindent
Here $\Cc_u,\, u = 1,\dots, l$, are non-intersecting oriented loops in $\C$.
The loops start at $+\)\infty$, go around 0, and return to $+\)\infty$.
For $u<v$, the loop $\Cc_u$ lies inside the loop $\Cc_v$. We fix a univalued
branch of the integrand by assuming that at the point of $\dl_l$ where all
numbers $-\)s_1, \dots,-\)s_l$ are positive we have $\arg\)(-\)s_u)=0$ for
$u=1,\dots,l$, and $\arg\)(s_u\]-s_v)=0$ for $1\le u<v\le l$.
\vsk.2>
The formula for $J_l(m)$ is well known:
\begin{equation}
\label{selberg}
J_l(m)\,=\,\ka^{\>l\)(l\)-1-\)m)/\ka}\;\prod_{j=0}^{l-1}\;
\frac{-2\pi i\;\Gm(1-1/\ka)}{\Gm(1+(m-j)/\ka)\,\Gm(1-(j+1)/\ka)}\;,
\end{equation}
for example, \cf. \cite{TV2}\), \cite{MTV}\).

\vsk.2>
It is easy to see that
\vvn.3>
\begin{align}
I_{a,b}(z;m_1,m_2,l_1,l_2)\, &{}=\,
e^{-bz/\ka}\,z^{-(2b^2-b\)(2l_2+m_2-m_1)+m_2l_2)/\ka}\,
e^{\pi i\)b\)m_1/\ka}\,\x{}
\notag
\\[4pt]
& \>{}\x\,(l_2\]-b)!\,b\)!\;C_{b}(m_1,m_2, l_1,l_2)\,
J_{l_2-b}(m_1)\,J_{b}(m_2)\,\bigl(\dl_{a,b}\>+\>\Oc(z^{-1})\bigr)\,.
\notag
\\[-10pt]
\notag
\end{align}
Substituting formulas for $C_{b}(m_1,m_2,l_1,l_2)$ and $J_{l}(m)$ we get the
lemma.
\end{proof}
\vsk.5>
\nt
{\bf Proof of Theorem~\ref{first}.}\enspace
By Lemma~\ref{lem:eq}, since ${B(m_1,m_2,l_1,l_2)\)=\)B(l_1,l_2,m_1,m_2)}$,
the matrices $\Ih(z;m_1,m_2,l_1,l_2)$ and $\Ih(z;l_1,l_2,m_1,m_2)$ obey the
same first order ordinary differential equation. Therefore,
\vvn.3>
$$
\Ih(z;m_1,m_2,l_1,l_2)\,=\,\Ih(z;l_1,l_2,m_1,m_2)\,X(l_1,l_2,m_1,m_2)\,.
\vv.3>
$$
Lemma~\ref{lem:as} implies that the connection matrix $X$ is simultaneously
lower and upper triangular with the unit diagonal. Hence $X$ is the identity.
Theorem~\ref{first} is proved.
\hfill\qed

\section{The $(\glt\>,\glt)$ duality for \KZ/ and dynamical differential
equations}
\subsection{The \KZ/ and dynamical differential equations for $\glt$}
Let $E_{i,j}$, $i,j=1,2$, be the standard generators of the Lie algebra
$\frak {gl}_2$ over $\C$.

We have the root decomposition \,$\glt=\n_+\oplus\h\oplus\n_-$ where
\vvn.4>
$$
\n_+\>=\,\C\cdot E_{1,2}\,,\qquad
\h\,=\,\C \cdot E_{1,1}\>\oplus\,\C \cdot E_{2,2}\,,\qquad
\n_-\>=\,\C \cdot E_{2,1}\,.
\vv.3>
$$
We identify the Cartan subalgebra $\h$ with $\C^2$ mapping
$\la=\la_1 E_{1,1}+\la_2 E_{2,2}$ to $(\la_1, \la_2)$.

\vsk.2>
The Casimir element $\Om\in\glt^{\ox 2}$ is defined by the formula
\vvn.5>
$$
\Omega\,=\,E_{1,1}\ox E_{1,1} + E_{2,2}\ox E_{2,2} +
E_{1,2}\ox E_{2,1} + E_{2,1}\ox E_{1,2}\,.
$$

\vsk.5>
Fix a non-zero complex number $\ka$. Let $V_1,\dots, V_n$ be \$\glt$-modules.
Set $V = V_1 \ox \dots \ox V_n$. The \em{\KZ/ operators}
$\nabla_a(z_1,\dots, z_n,\la_1,\la_2,\ka)$, $a=1,\dots, n$, acting on \V-valued
functions of complex variables $z_1,\dots, z_n, \la_1,\la_2$ are defined
by the formula
\vvn.6>
$$
\nabla_a(z_1,\dots,z_n,\la_1,\la_2,\ka)\,=\,
\ka\,\frac{\partial}{\partial z_a}\;-\>\sum_{\textstyle\fratop{b=1}{b \ne a}}^n
\frac{\Om^{(a,b)}}{z_a\] - z_b}\;-\,\la_1E_{1,1}^{(a)}-\la_2E_{2,2}^{(a)}\,,
\qquad a=1,\dots, n\,.
$$
The operator $\Om^{(a,b)} : V \to V$ acts as $\Om$ on $V_a\ox V_b$,
and as the identity on all other tensor factors. Similarly, the operator
$E_{i,i}^{(a)}$ acts as $E_{i,i}$ on $V_a$ and as the identity on all other
tensor factors. The \em{\KZ/ equations} for a function
$U(z_1,\dots, z_n,\la_1,\la_2)$ are
\vvn.6>
$$
\nabla_a(z_1,\dots,z_n,\la_1,\la_2,\ka)\,U(z_1,\dots, z_n,\la_1,\la_2)\,=\,0\,,
\qquad a=1,\dots, n\,.
$$
\vsk.5>
\ifUS\else\nt\fi
The \em{dynamical differential operators} $D_i(z_1,\dots,z_n,\la_1,\la_2,\ka)$,
$i=1,2$, acting on \V-valued functions of complex variables
$z_1,\dots, z_n, \la_1,\la_2$ are defined by the formulae
\ifUS\vvn.4>\fi
\begin{align}
D_1(z_1,\dots,z_n,\la_1,\la_2,\ka)\, &{}=\,
\ka\,\frac{\partial}{\partial\la_1}\,-\,
\frac1{\la_1\]-\la_2}\,(E_{2,1} E_{1,2}- E_{2,2})\,
-\>\sum_{a=1}^n\>z_a\>E_{1,1}^{(a)}\,,
\notag
\\[6pt]
D_2(z_1,\dots,z_n,\la_1,\la_2,\ka)\, &{}=\,
\ka\,\frac{\partial}{\partial\la_2}\,-\,
\frac1{\la_2\]-\la_1}\,(E_{2,1} E_{1,2}- E_{2,2})\,
-\>\sum_{a=1}^n\>z_a\>E_{2,2}^{(a)}\,.
\notag
\\[-10pt]
\notag
\end{align}
The \em {dynamical differential equations} for a function
$U(z_1,\dots, z_n,\la_1,\la_2)$ are
\vvn.5>
$$
D_i(z_1,\dots,z_n,\la_1,\la_2,\ka)\,U(z_1,\dots,z_n,\la_1,\la_2)\,=\,0\,,
\qquad i=1,2\,.
$$
\vsk.5>

\begin{theorem}
[\)\cite{FMTV}\)]\label{compatibility}
The \KZ/ equations and dynamical differential equations form a compatible
system of equations. Namely, we have
\vvn.5>
$$
[\nabla_a\,,\]\nabla_b]\,=\,0\,,\qquad [\nabla_a\>,D_i]\,=\,0\,\qquad
[D_1\>,D_2]\,=\,0\,,
\vv.4>
$$
for $a,b = 1,\dots, n$, and $i=1, 2$.
\end{theorem}
\goodbreak
\vsk.5>

The \KZ/ and dynamical differential operators preserve the weight decomposition
of $V$. Thus the \KZ/ and dynamical differential equations can be considered
as equations on a function $U(z_1,\dots,z_n,\alb\la_1,\la_2)$ taking values in
a given weight subspace of $V$.

\subsection{The duality}\label{S2}
For a complex number $m$, denote $M_m$ the \$\glt$-Verma module with highest
weight $(m,0)$ and highest weight vector $v_m$. The vectors $E_{2,1}^d\)v_m$,
$d\in\Zp$, form a basis in $M_m$.

For a non-negative integer $m$, denote $L_m$ the irreducible $\glt$-module
with highest weight $(m,0)$ and highest weight vector $v_m$. The vectors
$E_{2,1}^d\)v_m$, $d=0,\dots, m$, form a basis in $L_m$.

Let $m_1,l_1$ be complex numbers and $m_2, l_2$ non-negative integers
such that $m_1+m_2=l_1+l_2$. Consider the weight subspace
${(M_{m_1}\ox L_{m_2})}[l_1,l_2]$ of the tensor product
${M_{m_1}\ox L_{m_2}}$. The weight subspace has a basis
\vvn.3>
\begin{equation}
\label{Fbasis}
F^a(m_1,m_2,l_1,l_2)\,=\,
\frac1{(l_2-a)!\,a\)!}\;E_{2,1}^{l_2-a}\)v_{m_1}\ox E_{2,1}^a\)v_{m_2}\,,
\qquad a=0,\dots\min\)(m_2, l_2)\,.\!\!
\vv.4>
\end{equation}
There is a linear isomorphism
\vvn.3>
$$
\phi\):\)(M_{m_1}\ox L_{m_2})[l_1,l_2]\,\to\,
(M_{l_1} \ox L_{l_2})[m_1,m_2]\,,\ \quad
F^a(m_1,m_2,l_1,l_2)\,\mapsto\> F^a(l_1,l_2,m_1,m_2)\,.
$$
\vsk.5>

\begin{theorem}[\)\cite{TV1}\)]
Under the isomorphism $\phi$ the system of \KZ/ and dynamical differential
equations with values in $(M_{m_1}\ox L_{m_2})[l_1,l_2]$ turns into
the system of dynamical differential and \KZ/ equations with values
in $(M_{l_1}\ox L_{l_2})[m_1,m_2]${\rm:}
\vvn.5>
\begin{align}
\nabla_a(z_1,z_2,\la_1,\la_2,\ka)\;\phi\, &{}=\,
\phi\;D_{a}(\la_1,\la_2,z_1,z_2,\ka)\,,
\notag
\\[6pt]
D_a(z_1,z_2,\la_1,\la_2,\ka)\;\phi\, &{}=\,
\phi\;\nabla_a(\la_1,\la_2,z_1,z_2,\ka)\,,
\notag
\end{align}
for $a=1,2$.
\end{theorem}

\section{Hypergeometric solutions}
\subsection{Construction of hypergeometric solutions}
Define the \em{master} function by the formula
\vvn-.4>
\begin{align}
\Phi &{}(s_1,\dots,s_{l_2}, z_1,z_2,\la_1,\la_2; m_1, m_2)\,={}
\notag
\\[8pt]
& \>{}=\,e^{\la_1(m_1z_1+m_2z_2)\>-\>(\la_1-\la_2) \sum_{u=1}^{l_2}\!s_u}
(\la_1-\la_2)^{-l_2}\ (z_1-z_2)^{m_1m_2}\,\x{}
\notag
\\
& \kern3.56em {}\x\,\prod_{u=1}^{l_2}\,(z_1\]-s_u)^{-m_1}\>(z_2\]-s_u)^{-m_2}
\prod_{1\le u<v\le l_2}(s_u\]- s_v)^2\,.
\notag
\\[-12pt]
\notag
\end{align}
For $0\le a\le l_2$, define the \em{weight} function by the formula
\vvn.2>
$$
w_{l_2-a,\)a}(s_1,\dots,s_{l_2}, z_1,z_2)\,=\;
\Sym\biggl[\ \prod_{u=1}^{l_2-a}\,\frac 1{z_1-s_u}
\ \prod_{u=l_2-a+1}^{l_2}\,\frac 1{z_2-s_u}\ \biggr]\,.
\vv.2>
$$
Define an \$(M_{m_1}\ox L_{m_2})[l_1,l_2]$-valued differential \$l_2$-form
by the formula
\vvn.2>
$$
\om(s_1,\dots,s_{l_2},z_1,z_2)\,=\>\sum_{a=0}^{l_2} \,
w_{l_2-a,\)a}(s_1,\dots,s_{l_2}, z_1,z_2)\,ds^{l_2}\;
\frac 1{(l_2-a)!\,a\)!}\,E^{l_2-a}_{2,1}v_{m_1} \ox E^{a}_{2,1}v_{m_2}\ .
\vv.3>
$$

\vsk.2>
Fix a complex number $\ka$. Define an \$(M_{m_1}\ox L_{m_2})[l_1,l_2]$-valued
\vvn.5>
function by the formula
\begin{align}
\label{Udl}
U_\dl &{}(z_1,z_2,\la_1,\la_2)\,={}
\\[5pt]
&\>{}=\int_{\dl(z_1,z_2,\la_1,\la_2)}\!\!\!\!
\bigl(\Phi(s_1,\dots,s_{l_2},z_1,z_2, \la_1,\la_2; m_1, m_2)\bigr)^{\]1/\ka}\,
\om(s_1,\dots,s_{l_2},z_1,z_2)\,.
\notag
\end{align}
The function depends on the choice of integration chains
$\dl(z_1,z_2,\la_1,\la_2)$. We assume that for each $z_1,z_2,\la_1,\la_2$
the chain lies in $\C^{\>l_2}\}$ with coordinates $s_1,\dots,s_{l_2}$, and
the chains form a horizontal family of \$l_2$-dimensional homology classes
with respect to the multivalued function
$\Phi(s_1,\dots,s_{l_2}, z_1,z_2, \la_1,\la_2; m_1, m_2)^{1/\ka}\}$,
see a more precise statement below and in \cite{FMTV}.

\vsk.6>
\begin{theorem}
\label{th:hgeom}
For any choice of the horizontal family $\dl$, the function
$U_{\dl}(z_1,z_2,\la_1,\la_2)$ is a solution of the \KZ/ and dynamical
differential equations with values in $(M_{m_1}\ox L_{m_2})[l_1,l_2]$
and parameter $\ka$.
\end{theorem}
\vsk.1>
\nt
The theorem is a corollary of Theorem 3.1 in \cite{FMTV}.

\vsk.2>
The solutions of the \KZ/ and dynamical differential equations constructed
in this section are called the \em{hypergeometric solutions}.

\subsection{Proof of Lemma~\ref{lem:eq}}
Consider the ordinary differential equation for a function $\Psi(x)$ taking
values in $(M_{m_1}\ox L_{m_2})[l_1,l_2]$:
\vvn.4>
\begin{equation}
\label{additional}
\biggl(\,\ka\,\frac{\der}{\der x}\,-\,\frac{\Om - m_1m_2}x\,+\,
E^{(2)}_{2,2}\,\biggr)\,\Psi(x)\,=\,0\,.
\vv.3>
\end{equation}
One can directly check the following statement.

\begin{lemma}
\label{lem:UPsi}
Let $\Psi(x)$ be a solution of equation \Ref{additional}.
Then the function
\vvn.5>
\begin{align}
\label{UPsi}
U(z_1,z_2,\la_1,\la_2)\, &{}=\,
e^{\)(z_1\la_1(m_1-l_2) + z_1\la_2 l_2 + z_2\la_1 m_2)/\ka}\,\x{}
\\[4pt]
&\>{}\x\,(z_1-z_2)^{m_1m_2/\ka}\,(\la_1-\la_2)^{l_1l_2/\ka}\,
\Psi\bigl(\)-\)(\la_1\]-\la_2)\>(z_1\]-z_2)\bigr)
\notag
\\[-13pt]
\notag
\end{align}
is a solution of the \KZ/ and dynamical differential equations with values
in ${(M_{m_1}\}\ox\]L_{m_2})}[l_1,l_2]$ and parameter $\ka$. Moreover, formula
\Ref{UPsi} defines the one-to-one correspondence of solutions of the respective
equations.
\end{lemma}

For any $b=0,\dots\min\)(m_2,l_2)$ set
$$
\Ib_b(z;m_1,m_2,l_1,l_2)\,=\>\sum_{a=0}^{l_2}\,\,I_{a,b}(z; m_1,m_2,l_1,l_2)
\;\frac1{(l_2-a)!\,a\)!}\,E^{l_2-a}_{2,1}\)v_{m_1} \ox E^{a}_{2,1}\)v_{m_2}\,.
$$
Here the actual range of summation is until $a=\min\)(m_2,l_2)$,
since $E^{a}_{2,1}\)v_{m_2}=0$ for $a>m_2$.

\begin{lemma}
The function
\vvn.2>
\begin{align}
U_b&{}(z_1,z_2,\la_1,\la_2)\,=\,
e^{\)(z_1\la_1(m_1-l_2) + z_1\la_2 l_2 + z_2\la_1 m_2)/\ka}\,\x{}
\label{Ub}
\\[4pt]
&\,{}\x\,(z_1-z_2)^{m_1m_2/\ka}\,(\la_1-\la_2)^{l_1l_2/\ka}\,
\Ib_b\bigl(\)-\)(\la_1\]-\la_2)\>(z_1\]-z_2);m_1,m_2,l_1,l_2)\bigr)
\notag
\\[-13pt]
\notag
\end{align}
is a solution of the \KZ/ and dynamical differential equations with values
in ${(M_{m_1}\}\ox\]L_{m_2})[l_1,l_2]}$ and parameter $\ka$.
\end{lemma}
\begin{proof}
Change the variables $t_1,\dots t_{l_2}$ in the integrals
\ifUS\vvn.5>\fi
$I_{a,b}(z;m_1,m_2,l_1,l_2)$ to the new variables
$$
s_u\,=\,t_u\)(\la_1-\la_2)^{-1}+\)z_1\,,\qquad u=1,\dots,l_2\,.
\vv.5>
$$
Then the right hand side of \Ref{Ub} takes the form of integral \Ref{Udl}
with the chain $\dl(z_1,z_2,\alb \la_1,\la_2)$ induced by the chain
$\gm_b\bigl(\)-\)(\la_1\]-\la_2)\>(z_1\]-z_2)\bigr)$ under the described
change of integration variables. One can see that $\dl$ is a horizontal family
of chains. Therefore, the statement follows from Theorem~\ref{th:hgeom}.
\end{proof}

To complete the proof of Lemma~\ref{lem:eq} observe that equation
\Ref{additional} written in basis \Ref{Fbasis} of the weight subspace
$(M_{m_1}\ox L_{m_2})[l_1,l_2]$ becomes equation \Ref{eq:Ihat}\). Hence,
Lemma~\ref{lem:UPsi} yields Lemma~\ref{lem:eq}.

\section{Remarks on asymptotics with respect to dimension}

Let $\ka$ be a positive number. Consider a 4-tuple $m_1, m_2, l_1, l_2$
of complex numbers such that $m_2, l_2$ are natural numbers and
$m_1+m_2 = l_1+l_2$. For any integers $a, b$, admissible with respect
to the pair $m_2, l_2$, introduce an $l_2$-dimensional integral
$K_{a,b}(z;m_1,m_2,l_1,l_2)$ by the formula
\vvn.1>
$$
K_{a,b}(z; m_1, m_2, l_1, l_2)\,=
\int_{\gm_{l_2,\)b}(z)\!}\!\Phi_{l_2}(t, z; m_1, m_2)^{1/\ka}
\,w_{l_2-a,\)a}(t, z)\,dt^{\>l_2}\,,
\vv.2>
$$
where the cycle $\gm_{l_2,b}(z)$ and the functions $\Phi_{l_2}(t,z;m_1,m_2)$,
$w_{l_2-a,\)a}(t, z)$ are defined in the Introduction. Recall that
$\gm_{l_2,b}(z)$ is the \$l_2$-dimensional multi-loop with $b$ loops
going around $z$ and the remaining $l_2-b$ loops going around 0.

\begin{corollary}
We have
$$
K_{a,b}(z, m_1,m_2,l_1,l_2)\,=\,e^{\pi i\)b\)(l_2-m_2)}\;
\frac{(l_2\]-b)!\,J_{l_2-b}(m_1)\>J_b(m_2)}
{(m_2\]-b)!\,J_{m_2-b}(l_1)\>J_b(l_2)}\;K_{a,b}(z, l_1,l_2,m_1,m_2)\,.
$$
\end{corollary}
\vsk.6>

Assume that the numbers $m_1, m_2$ are fixed and $l_2$ tends to infinity.
Then the integral $K_{a,b}(z, m_1,m_2,l_1,l_2)$ is a hypergeometric integral
of growing dimension $l_2$ with fixed exponents $m_1, m_2$. In the right hand
side we have $J_{l_2-b}(m_1)$, the Selberg integral of growing dimension;
$J_{m_2-b}(l_1), J_b(l_2)$, the Selberg integrals of fixed dimension with
growing exponents $l_1, l_2$; $K_{a,b}(z, l_1,l_2,m_1,m_2)$, the hypergeometric
integral of fixed dimension $m_2$ with growing exponents $l_1, l_2$.

\vsk.1>
The Selberg integrals are given explicitly by formula \Ref{selberg}.
The asymptotics of the integral $K_{a,b}(z, l_1,l_2,m_1,m_2)$
can be calculated by the steepest descent method.

Here is an example. If $m_2=1$, then the integral $K_{a,b}(z,l_1,l_2,m_1,m_2)$
is one-dimen\-sional and its asymptotics are described as follows.

Consider the integral
\vvn.1>
$$
I_C\,=\,\int_{\!C\,} e^{-t}\,(-t)^{M+a}\,(z-t)^{-M}\,dt
\vv.3>
$$
where $a\in\C$, $M$ is a real parameter, and $C$ is one of two loops $C',\>C''$
shown in the picture.

\ifUS
\vsk-.2>
\else
\vsk2>\goodbreak\vsk-2.8>\fi
\vbox{
\begin{center}
\begin{picture}(224,85)
\put(30,35){\circle*{3}}
\put(37,30){$z$}
\put(40,15){\circle*{3}}
\put(47,10){$0$}
\put(40,15){\oval(16,16)[l]}
\put(40,7){\vector(1,0){160}}
\put(40,23){\line(1,0){160}}
\put(205,20){$C''$}
\qbezier(40,0)(10,0)(10,32)
\qbezier(10,32)(10,65)(40,65)
\put(40,0){\vector(1,0){160}}
\put(40,65){\line(1,0){160}}
\put(205,62){$C'$}
\end{picture}
\vsk>
Picture 3. The contours $C',\>C''$.
\end{center}}

\vsk.7>
\noindent
We assume that at the point of $C$ where $t$ is negative we have
$\arg\)(-\)t)=0$ and $\arg\)(z-t)\in(0,\pi)$.
We are interested in asymptotics of the integral as $M$ tends to $+\)\infty$.

\vsk.5>
\begin{lemma} Assume that $\Im z>0$ and $\arg\)(-\)z)\in (-\)\pi, 0)$.
Assume that $\arg M=0$ and $M$ tends to $+\)\infty$. Then
\vvn.5>
\begin{gather}
I_{C'}\>=\,i\>\pi^{1/2}\)(-\)zM)^{(2a+1)/4}\)
\exp\)\bigl(\)2\)(-\)zM)^{1/2}\]-z/2\bigr)\)\bigl(1+\Oc(M^{-1/2})\bigr)\,,
\notag
\\[8pt]
I_{C''}\>=\,\pi^{1/2}\)(e^{2\pi i\)(M+\)a)}-1)\,
(-\)zM)^{(2a+1)/4}\)\exp\)\bigl(-2\)(-\)zM)^{1/2}\]-z/2-\pi i\)a\bigr)\)
\bigl(1+\Oc(M^{-1/2}) \bigr)\,.\ifUS\else\kern-7pt\fi
\notag
\end{gather}
\end{lemma}
\vsk.5>

\nt
To prove the lemma one changes the integration variable, $u=t\)M^{-1/2}$, then
\vvn.8>
$$
I_C\,=\,M^{\)(a+1)/2}\)\int_{\!C\,}
\exp\)\bigl(M^{\)1/2}(-u + z/u)+ z^2\]/2u^2\bigr)\,(-u)^a\,du
\;\bigl(1+\Oc(M^{-1/2})\bigr)\ ,
\vv.2>
$$
and the integral becomes a standard integral of the steepest descent method,
whose asymptotics come from critical points of the function $-\)u+z/u$.

\vsk>


\begin{thebibliography}{MTV}
\frenchspacing
\small

\bibitem[FMTV]{FMTV}
Felder, G., Markov, Y., Tarasov, V., Varchenko, A.: Differential equations
compatible with \KZ/ equations, Math. Phys. Anal. Geom. {\bf 3} (2000), no. 2,
139--177.

\bibitem[MTV]{MTV}
Markov, Y., Tarasov, V., Varchenko, A.: The determinant of a hypergeometric
period matrix, Houston J. of Math. {\bf 24} (1998) no. 2, 197--220.

\bibitem[TV1]{TV1}
Tarasov, V., Varchenko, A.: Duality for Knizhnik-Zamolodchikov and
dynamical equations, 
Acta Appl. Math. {\bf 73} (2002), no. 1-2, 141--154.

\bibitem[TV2]{TV2}
Tarasov, V., Varchenko, A.: Selberg type integrals associated with $sl_3$,
Preprint (2003) 1--11, {\tt math.QA/0302148}.

\end{thebibliography}
\end{document}